\documentstyle[12pt]{article}
\begin{document}

\def\dsp{\displaystyle}
\def\Rr{{bf R}} 
\def\Zz{bf Z}
\def\Nn{bf N}
\def\get{\hbox{{\goth g}$^*$}}
\def\g{\gamma}
\def\om{\omega}
\def\r{\rho}
\def\a{\alpha}
\def\s{\sigma}
\def\vfi{\varphi}
\def\l{\lambda}
\def\implique{\Rightarrow}
\def\o{{\circ}}
\def\Diff{\hbox{\rm Diff}}
\def\S1{\hbox{\rm S$^1$}}
\def\Hom{\hbox{\rm Hom}}
\def\Vect{\hbox{\rm Vect}}
\def\const{\hbox{\rm const}}
\def\ad{\hbox{\hbox{\rm ad}}}
\def\semid{\hbox{\bb o}}
\def\blanc{\hbox{\ \ }}

\def\pds#1,#2{\langle #1\mid #2\rangle} 
\def\f#1,#2,#3{#1\colon#2\to#3} 

\def\hfl#1{{\buildrel{#1}\over{\hbox to
12mm{\rightarrowfill}}}}

\title{Space of linear differential operators on the real line as a module over
the Lie algebra of vector fields}

\author{H. Gargoubi, V.Yu. Ovsienko\\
{\small CNRS, Centre de Physique Th\'eorique}
\thanks{CPT-CNRS, Luminy Case 907,
F--13288 Marseille, Cedex 9
FRANCE
}
}

\date{}

\maketitle

{\abstract{Let ${\cal D}^k$ be the space of $k$-th order linear differential
operators on ${\bf R}$: $A=a_k(x)\frac{d^k}{dx^k}+\cdots+a_0(x)$.
We study a natural 1-parameter family of $\Diff(\bf R)$- (and $\Vect(\bf R)$)-modules 
on ${\cal D}^k$.
(To define this family, one considers arguments of differential
operators as tensor-densities of degree $\lambda$.)
In this paper we solve the problem of isomorphism 
between $\Diff(\bf R)$-module structures on ${\cal D}^k$
corresponding to different values 
of  $\lambda$. 
The result is as follows: for $k=3$ $\Diff(\bf R)$-module 
structures on ${\cal D}^3$
are isomorphic to each other
for every values of
$\lambda\not=0,\;1,\;{1\over 2},\;{1\over 2}\pm \frac{\sqrt 21}{6}$,
in this case
there exists a unique (up to a constant)
intertwining operator $T:{\cal D}^3\to{\cal D}^3$.
In the higher order case $(k\geq 4)$ $\Diff(\bf R)$-module structures on ${\cal D}^k$
corresponding to two different values of the degree: $\lambda$ and $\lambda^{\prime}$,
are isomorphic if and only if $\lambda+\lambda^{\prime}=1$.
}}                                                               																																																																									

\thispagestyle{empty}

\vfill\eject

\section{Introduction}
Space of linear differential operators on a manifold $M$ has 
various algebraic structures: the structure of associative algebra and of Lie algebra,
in the 1-dimensional case it can be considered as an
infinite-dimensional Poisson space (with respect to so-called Adler-Gelfand-Dickey
bracket).

\vskip 0,3cm

{\bf 1.1 $\Diff(M)$-module structures}.
One of the basic structures on the
space of linear differential operators
is a natural family of module structures over the
group of diffeomorphisms $\Diff(M)$ (and
of the Lie  algebra of vector fields $\Vect(M)$).
These $\Diff(M)$- (and $\Vect(M)$)-module structures
are defined if one considers the arguments
of differential operators as {\it tensor-densities}
of degree $\lambda$ on $M$.

\vskip 0,3cm

In this paper we consider the space of differential operators on 
${\bf R}$.\footnote{Particular cases of actions of $\Diff(\bf R)$ and $\Vect(\bf R)$
on this space
were considered by classics (see \cite{car}, \cite{wil}).
The well-known example is the Sturm-Liouville operator
$\frac{d^2}{dx^2}+a(x)$ acting on $-\frac{1}{2}$-densities
(see e.g. \cite{car}, \cite{wil}, \cite{tur}). 
 Already
this simplest case
leads to interesting geometric structures
and is related to so-called Bott-Virasoro group (cf. \cite{kir}, \cite{seg}).}
Denote ${\cal D}^k$ the space of $k$-th order linear differential 
operators:
\begin{equation}
A(\phi)=a_k(x)\frac{d^k\phi}{dx^k}+\cdots+a_0(x)\phi
\label{ope}
\end{equation}
where $a_i(x),\phi(x)\in C^{\infty}(\bf R)$.

Define a 1-parameter family of $\Diff(\bf R)$-module structures on 
$C^{\infty}({\bf R})$ by:
$$
g^*_{\lambda}\phi:=\phi\circ g^{-1}\cdot\left(\frac{dg^{-1}}{dx}\right)^{-\lambda}
$$
where $\lambda \in \bf R$ (or $\lambda \in \bf C$) is a parameter.
Geometrically speaking, $\phi$ 
is a {\it tensor-density} of degree $-\lambda$ :
$$
\phi =\phi (x)(dx)^{-\lambda }.
$$

A 1-parameter family of actions of  $\Diff(\bf R)$ on the space of differential
operators (\ref{ope}) is defined by:
$$
g(A)=g^*_{\lambda}A(g^*_{\lambda})^{-1}.
$$

Denote ${\cal D}^k_{\lambda}$ 
the space of operators (\ref{ope})
endowed with the defined $\Diff(\bf R)$-module structure.

Infinitesimal version of this action defines a 1-parameter family of 
$\Vect(\bf R)$-module structures
on ${\cal D}^k$ (see Sec. 3 for details).

\vskip 0,3cm

{\bf 1.2 The problem of isomorphism}.
Let $M$ be a manifold, $\dim M\geq 2$.
The problem of isomorphism of
$\Diff(M)$- (and $\Vect(M)$)-module structures
for different values of $\lambda$
was stated  in \cite{duv}
and saved in the case of second order differential operators.
In this case,
different  $\Diff(M)$-module structures are isomorphic to each other
for every $\lambda$ except 3 critical values: $\lambda=0,-{1\over 2},-1$
(corresponding to differential operators on: functions, $\frac{1}{2}$-densities
and volume forms respectively).

Geometric quantization gives an example of such a special $\Diff(M)$-module:
 differential operators are considered as acting on $\frac{1}{2}$-densities
(see \cite{kos}).

Recently P.B.A. Lecomte, P. Mathonet, and E. Tousset 
\cite{lec} showed that in the case of differential operators of order $\geq 3$,
$\Diff(M)$-modules
corresponding to $\lambda$ and $\lambda^{\prime}$-densities
are isomorphic if and only if $\lambda+\lambda^{\prime}=1$.
The unique isomorphism in this case
is given by conjugation
of differential operators.

These results solve the
problem of isomorphism in the multi-dimensional case.

It was shown in \cite{duv}, \cite{lec}, that the case 
$\dim M =1$ ($M={\bf R}$ or $S^1$) is 
particular.
It is reacher in algebraic structures and
therefore is of a special interest.

In this paper we solve the problem of isomorphism
of $\Diff(\bf R)$-modules ${\cal D}^k_{\lambda}$
for any $k$. The result is as follows.

1) The modules
${\cal D}^3_{\lambda}$ of {\it third order}
differential operators (\ref{ope})
are isomorphic to each other for all values of $\lambda$ except
5 critical values:
$$
\{0,\;\;-1,\;\;-{1\over 2},\;\;\;-{1\over 2}+ \frac{\sqrt 21}{6},\;\;\;-{1\over 2}- \frac{\sqrt 21}{6}\}.
$$
(this result was announced in \cite{duv}).

2) The $\Diff({\bf R})$-modules ${\cal D}^k_{\lambda}$ and ${\cal D}^k_{\lambda^{\prime}}$
on the space of differential operators (\ref{ope}) of order
 $k\geq 4$ are isomorphic if and only if
$\lambda+\lambda^{\prime}=-1$.

\vskip 0,3cm

{\bf 1.3 Intertwining operator}.
The most important result of the paper
is a construction of the unique (up to a constant) equivariant linear operator
on the space of third order differential operators:
\begin{equation}
T:{\cal D}^3_{\lambda}\to{\cal D}^3_{\mu}
\label{T}
\end{equation}
for $\lambda,\mu\not=0,-1,-{1\over 2},\;-{1\over 2}\pm \frac{\sqrt 21}{6}$,
see the explicit formul{\ae} (\ref{int}), (\ref{TTT}) and (\ref{exp}) below.
It has nice geometric and algebraic properties and
seems to be an interesting object to study.

Operator $T$ is an analogue
of the {\it second order Lie derivative} from \cite{duv}
intertwining different $\Diff(M)$-actions on the space of second order differential operators
on a multi-dimensional manifold $M$.

\vskip 0,3cm

{\bf 1.3 Normal symbols}.
The main tool of this paper is the notion of a {\it normal symbol},
which we define in the case of 4-th order differential operators.
We define a $sl_2$-equivariant way to associate a polynomial function 
of degree 4 on $T^*\bf R$
to a differential
operator $A\in {\cal D}^4_{\lambda}$.
In the case of second order operators the notion of
normal symbol was defined in \cite{duv}.
This construction is related with the results of
\cite{cmz}. We discuss the geometric properties of the normal symbol
and its relations to
the intertwining operator (\ref{T}).

\section{Main results}

We formulate here the main results of this paper,
all the proofs will be given in Sec. 3-7. 

\vskip 0,3cm

{\bf 2.1 Classification of $\Diff(\bf R)$-modules}.
First, remark that for each value of $k$,
there exists an isomorphism of $\Diff(\bf R)$-modules:
$$
{\cal D}^k_{\lambda} \cong {\cal  D}^k_{-1-\lambda}.
$$
It is given by {\it conjugation} $A\mapsto A^*$:
$$
A^*
=
\sum_{i=1}^{k}(-1)^i\frac{d^i}{dx^i}\circ a_i(x)
$$

The following two theorems give a solution of the 
problem of isomorphism 
of $\Diff(\bf R)$-modules ${\cal D}^k_{\lambda}$
on space ${\cal D}^k$.

The first result 
was announced on \cite{duv}:

\proclaim Theorem 1.
(i) All the $\Diff(\bf R)$-modules ${\cal D}^3_{\lambda}$ 
with $\lambda\not=0,-1,-\frac{1}{2},-\frac{1}{2}+ \frac{\sqrt 21}{6},
-\frac{1}{2}- \frac{\sqrt 21}{6}$ 
are isomorphic to each other. 
\hfill\break
(ii) The modules 
$
{\cal D}^3_0,\;{\cal D}^3_{-\frac{1}{2}},\;{\cal D}^3_{-\frac{1}{2}+\frac{\sqrt 21}{6}}
$
are not isomorphic to ${\cal D}^3_{\lambda}$ for general $\lambda$.\par

It follows from the general isomorphism $*:{\cal D}^k_{\lambda}\cong{\cal D}^k_{-1-\lambda}$,
that 
$$
{\cal D}^3_0\cong{\cal D}^3_{-1}\;\;\;\hbox{and}\;\;\;
{\cal D}^3_{-\frac{1}{2}+\frac{\sqrt 21}{6}}\cong
{\cal D}^3_{-\frac{1}{2}-\frac{\sqrt 21}{6}}.
$$
Therefore, there exist 4 non-isomorphic 
$\Diff(\bf R)$-module structures on the space ${\cal D}^3$.

\proclaim Theorem 2. For $k\geq 4$, the $\Diff(\bf R)$-modules  
${\cal D}^k_{\lambda}$ and ${\cal D}^k_{\lambda^{\prime}}$
are isomorphic if and only if $\lambda+\lambda^{\prime}=-1$.\par

This result shows that operators of order 3 play a special role
in the 1-dimensional case (as operators of order 2 in the case of a 
manifold of dimension $\geq 2$, cf. \cite{duv}, \cite{lec}). 

\vskip 0,3cm

{\bf 2.2 Intertwining operator $T$}.

\proclaim Theorem 3. For 
$\lambda,\mu\not=0,-1,-\frac{1}{2},-\frac{1}{2}\pm \frac{\sqrt 21}{6}$
there
exists a unique (up to a constant)
 isomorphism 
of $\Diff(\bf R)$-modules 
${\cal D}^3_{\lambda}$ and ${\cal D}^3_{\mu}$.\par

Let us give an explicit formula for the operator (\ref{T}).

Every differential operator of order 3
can be written (not in a canonical way) as a linear combination of
four operators: 
\hfill\break
1) a zero order operator of multiplication by a function: 
$\phi(x)\mapsto\phi(x)f(x)$,
\hfill\break
2) a first order operator of Lie derivative:
$$
L_X^{\lambda}=X(x)\frac{d}{dx}-\lambda X^{\prime}(x),
$$
where $ X^{\prime}=\frac{dX}{dx}$,
\hfill\break
3) symmetric ``anti-commutator'' of Lie derivatives:
$$
[L_X^{\lambda},L_Y^{\lambda}]_+:=
L_X^{\lambda}\circ L_Y^{\lambda}+
L_Y^{\lambda}\circ L_X^{\lambda}
$$
\hfill\break
4) symmetric third order expression:
$$
 [L_X^{\lambda},L_Y^{\lambda},L_Z^{\lambda}]_+:=
Sym_{X,Y,Z}
(L_X^{\lambda}\circ L_Y^{\lambda}\circ L_Z^{\lambda})
$$
for some vector fields $X(x)\frac{d}{dx},Y(x)\frac{d}{dx},Z(x)\frac{d}{dx}$.

\proclaim Theorem 4. The following formula:
\begin{equation}
\matrix{
T(f)=
\displaystyle\frac{\mu(\mu+1)(2\mu+1)}{\lambda(\lambda+1)(2\lambda+1)}f\hfill\cr\noalign{\bigskip}
T(L_X^{\lambda})=
\displaystyle\frac{3\mu^2+3\mu-1}{3\lambda^2+3\lambda-1}L_X^{\mu}\hfill\cr\noalign{\bigskip}
T([L_X^{\lambda},L_Y^{\lambda}]_+)=
\displaystyle\frac{2\mu+1}{2\lambda+1}[L_X^{\mu},L_Y^{\mu}]_+\hfill\cr\noalign{\bigskip}
T([L_X^{\lambda},L_Y^{\lambda},L_Z^{\lambda}]_+)=
[L_X^{\mu},L_Y^{\mu},L_Z^{\mu}]_+\hfill\cr
}
\label{int}
\end{equation}
defines an intertwining operator (\ref{T}). \par

A remarkable fact is that
the formula (\ref{int}) does not depend on the choice of
$X,Y,Z$ and $f$ representing the third order operator.
(Indeed,
the formul{\ae} (\ref{TTT}) and (\ref{exp}) below give the expression of $T$
directly
in terms of coefficients of differential operators.)
Moreover, this property fixes the coefficients in (\ref{int})
in a unique way (up to a constant).

\vskip 0,3cm
 
{\bf Remarks}. 1) In the case of multi-dimensional manifold $M$,
almost all $\Diff(M)$-module structures 
on the space of second order differential operators
are isomorphic to each other and the corresponding isomorphism
is unique (up to a constant) \cite{duv};
there is no isomorphism between different
$\Diff(M)$-module structures on the space of third order
operators, except the conjugation \cite{lec}.

2) The formula (\ref{int}) gives an idea that
it would be interesting to study
the commutative algebra structure 
(defined by the anti-commutator)
on the Lie algebra of all differential operators.

\section {Action of $\Vect(\bf R)$ on space ${\cal D}^4$}  

To prove Theorems 1-4, it is sufficient to consider
only the $\Vect(\bf R)$-action on ${\cal D}^k$. Indeed, since
the $\Diff(\bf R)$-action on the space of differential
operators is local, therefore,
the two properties: of $\Vect(\bf R)$- and of $\Diff(\bf R)$-equivariance 
are equivalent.

\vskip 0,3cm
 
{\bf 3.1 Definition of the family of $\Vect(\bf R)$-actions}.
Let $\Vect(\bf R)$ be the Lie algebra of smooth vector fields on $\bf R$:
$$
X=X(x)\frac{d}{dx}
$$
with the commutator
$$
[X(x)\frac{d}{dx},Y(x)\frac{d}{dx}]=
(X(x)Y^{\prime}(x)-X^{\prime}(x)Y(x))\frac{d}{dx},
$$
where $X^{\prime}=dX/dx$.

The action of $\Vect(\bf R)$ on space ${\cal D}^k$ is defined by:   
$$
\ad L_{X}^{\lambda}(A):=L_{X}^{\lambda}\circ A-A\circ L_{X}^{\lambda}
$$
where
$$
L_X^{\lambda}\phi =X(x)\phi^{\prime}(x)-\lambda X^{\prime}(x)\phi(x)
$$
The last formula defines a 1-parameter family of 
$\Vect(\bf R)$-actions on $C^{\infty}(\bf R)$.

One obtains a 1-parameter family of 
$\Vect(\bf R)$-modules on ${\cal D}^k$.

\vskip 0,3cm

{\bf  Notation}. 
1. The operator $L_X^{\lambda}$ is 
called the {\it operator of Lie derivative} of tensor-densities of degree $-\lambda$.
Denote ${\cal F}_{\lambda}$ the corresponding $\Vect(\bf R)$-module structure
on $C^{\infty}(\bf R)$.

2. As in the case of $\Diff({\bf R})$-module structures, we
denote ${\cal D}^k_{\lambda}$
space ${\cal D}^k$
as a $\Vect(\bf R)$-module.

 \vskip 0,3cm

{\bf 3.2 Explicit formula}.
Let us calculate explicitly the action of
 Lie algebra $\Vect(\bf R)$ on space ${\cal D}^4$. 
Given a differential operator $A\in{\cal D}^4$,
let us use the following notation for the $\Vect(\bf R)$-action
$\ad L_X$:
$$
\ad L_X(A)=a_4^X(x)\frac{d^4}{dx^4}+a_3^X(x)\frac{d^3}{dx^3}+a_2^X(x)\frac{d^2}{dx^2}+
a_1^X(x)\frac{d}{dx}+a_0^X(x).
$$

\proclaim Lemma 3.1. The action $\ad L_{X}^{\lambda}$ of $\Vect(\bf R)$  
on space ${\cal D}^4$ is given by : 
\begin{equation}
\matrix{ 
a_4^X
=
L_X^4(a_4) \hfill \cr\noalign{\smallskip} 
a_3^X
=
L_X^3(a_3) + 2(2\lambda -3)a_4X'' \hfill  \cr\noalign{\smallskip}
a_2^X
=
L_X^2(a_2) + 3(\lambda -1)a_3X'' + 2(3\lambda -2)a_4X''' \hfill \cr\noalign{\smallskip}
a_1^X
=
L_X^1(a_1) + (2\lambda -1)a_2X'' + (3\lambda -1)a_3X''' + (4\lambda -1)a_4X^{IV} \hfill  \cr\noalign{\smallskip}
a_0^X  
=
L_X^{0}(a_0) + \lambda (a_1X'' + a_2X''' + a_1X^{IV} + a_0X^{V}) \hfill \cr\noalign{\smallskip}
}
\label{coe} 
\end{equation}\par

{\bf Proof}. One gets easily the formula (\ref{coe}) from the definition: 
$$ 
\matrix{
\ad L_X^{\lambda}(A) = [L_X^{\lambda},A] =&
 (X\displaystyle\frac{d}{dx} -\lambda X') (a_4\frac{d^4}{dx^4}+\cdots+a_0)
\cr\noalign{\smallskip}
&-(a_4\displaystyle\frac{d^4}{dx^4}+\cdots+a_0)
(X\frac{d}{dx} -\lambda X')\hfill\cr\noalign{\smallskip} 
}
$$ 

\vskip 0,3cm

{\bf 3.3 Remarks}. It is convenient to interpret the action (\ref{coe}) as
 a {\it deformated} standard action 
  of $\Vect(\bf R)$ on the direct sum:
$$
{\cal F}_4\oplus{\cal F}_3\oplus{\cal F}_2
  \oplus{\cal F}_1\oplus {\cal F}_0.
$$ 
(given by the first term of the right hand side of each equality
in the formula (\ref{coe})). This interpretation 
is the motivation of the main construction of Sec. 4, it will be
discussed in Sec. 7.2.

  The main idea of proof of Theorems 1 and 2 is
 to find some {\it normal form} (cf. \cite{duv}) of the coefficients 
$a_4(x),\dots,a_0(x)$
  for 4-order differential operators on $\bf R$ which reduce the 
  action (\ref{coe}) to a canonical form.

 \section{Normal form of a symbol}

   It is convenient to represent differential operators as
   polynomials on the cotangent bundle.
 The standard way to define a (total) symbol
of an an operator (\ref{ope})
is to associate to 
   $A$ the polynomial
$$
P_A(x,\xi) =\sum_{i=0}^{k}\xi^ia_i(x),
$$
 
on $T^*{\bf R}\cong{\bf R}^2$  (where $\xi$ is a coordinate on the fiber). 
   However, this formula depends on coordinates, only the higher term 
   $\xi^ka_k(x)$ of $P_A$ (the principal symbol) has a geometric sense.

\vskip 0,3cm
   
{\bf 4.1 The main idea}.
Lie algebra $\Vect(\bf R)$ naturally acts 
on $C^{\infty}(T^*{\bf R})$ (it acts on the cotangent bundle).

 Consider a linear differential operator
$A\in {\cal D}^4$. 
Let us look for a natural definition of a symbol of $A$ in the following form:
$$
   \overline{P} _A(x,\xi)=\xi^4\bar a_4(x) +\xi^3\bar a_3(x)+
\xi^2\bar a_2(x)  +\xi\bar a_1(x)+\bar a_0(x),
$$  
where the functions $\bar a_i(x)$ are linear expressions of the coefficients 
  $a_i(x)$ and their derivatives.

Any symbol $\overline{P}(x,\xi)$ can be considered as a linear mapping
$$
{\cal D}^4\to{\cal F}_4\oplus{\cal F}_3\oplus{\cal F}_2
  \oplus{\cal F}_1\oplus{\cal F}_0.
$$
Indeed, the
Lie algebra $\Vect(\bf R)$ acts on each coefficient $\bar a_i(x)$ of the polynomial 
$\overline{P} _A(x,\xi)$ as on a tensor-density of degree $-i$:
$$
L_X(\overline{P} _A)=\sum_{i=0}^{4}\xi^iL^i_X(\bar a_i).
$$
However, there is no such a mapping which is $\Vect(\bf R)$-equivariant.

\proclaim 4.2 Definition.  The {\it normal symbol} of $A\in {\cal D}^4_{\lambda}$ as a polynomial
$\overline{P} _A(x,\xi)$
 such that the linear mapping
$
A\mapsto\overline{P} _A
$
is equivariant with respect to the subalgebra $sl_2\subset\Vect(\bf R)$
generated by the vector fields
$$
\{\frac{d}{dx} ,\;\; x\frac{d}{dx} , \;\; x^2\frac{d}{dx}\}.
$$\par

  \proclaim Proposition. 4.1.  (i) The following formula defines a normal 
  symbol of a differential operator $A\in {\cal D}^4_{\lambda}$:
  \begin{equation}  
\matrix{  
\bar a_4
=&
a_4  \hfill \cr\noalign{\smallskip}
\bar a_3 
=&
a_3 + {1\over 2}(2\lambda -3)a_4' \hfill \cr\noalign{\smallskip}
\bar a_2 
=& 
a_2 + (\lambda -1)a_3' + {2\over 7}(\lambda -1)(2\lambda -3)a_4'' \hfill \cr\noalign{\smallskip}
\bar a_1  
=& 
a_1 + {1\over 2}(2\lambda -1)a_2' + {3\over 10}(\lambda -1)(2\lambda 
-1)a_3''  \hfill \cr\noalign{\smallskip}
&
\;\;\;\;+ {1\over 15}(\lambda -1)(2\lambda -1)(2\lambda -3)a_4'''     \hfill \cr\noalign{\smallskip}
\bar a_0 
=& 
a_0 + \lambda a_1' + {1\over 3}{\lambda}(2\lambda -1)a_2'' + 
{1\over 6}{\lambda}(\lambda -1)(2\lambda -1)a_3'''    \hfill \cr
&
\;\;\;\;+
{1\over 30}{\lambda}(\lambda -1)(2\lambda -1)(2\lambda -3)a_4^{(IV)}\hfill \cr
}
\label{cum} 
\end{equation}
\hfill\break
(ii) The normal symbol is defined uniquely (up to multiplication
of each function $\bar a_i(x)$ by a constant).\par

{\bf Proof}.  Direct calculation shows that the 
$\Vect(\bf R)$-action $ad L^{\lambda}$ on ${\cal D}^4$ given by the 
formula (\ref{coe}) reads in terms of $\bar a_i$ as: 
\begin{equation}
\matrix{ 
\bar a_4 ^X
=&
L_X^{4}(\bar a_4) \hfill \cr\noalign{\smallskip} 
\bar a_3 ^X 
=& 
L_X^{3}(\bar a_3)  \hfill  \cr\noalign{\smallskip}
\bar a_2 ^X
=& 
L_X^{2}(\bar a_2) + {2\over 7}(6{\lambda}^2+6\lambda -5)J_3(X,\bar a_4) \hfill \cr\noalign{\smallskip}
\bar a_1 ^X
=& 
L_X^{1}(\bar a_1) + {2\over 5}(3{\lambda}^2+3\lambda -1)J_3(X,\bar a_3)  \hfill  \cr\noalign{\smallskip}
&
\;\;\;\;\;\;\;\;\;\;\;\;
+ {1\over 6}{\lambda}(\lambda +1)(2\lambda +1)J_4(X,\bar a_4)    \hfill  \cr\noalign{\smallskip}
\bar a_0 ^X 
=& 
L_X^{0}(\bar a_0)   
+ {2\over 3}{\lambda}(\lambda +1)J_3(X,\bar a_2)  \hfill  \cr\noalign{\smallskip}
&
\;\;\;\;\;\;\;\;\;\;\;\; 
 + {1\over 6}{\lambda}(\lambda +1)(2\lambda +1)J_4(X,\bar a_3)    \hfill  \cr\noalign{\smallskip}
&
\;\;\;\;\;\;\;\;\;\;\;\;
+ {1\over 420}{\lambda}(\lambda +1)(12{\lambda}^2+12\lambda +11)J_5(X,\bar a_4)  \hfill 
}
\label{nac} 
\end{equation}
where $\bar a_i ^X$ are coefficients of the normal symbol of the operator
$adL_X^{\lambda}(A)$ and the expressions $J_m$ are:
$$
\begin{array}{l}  
J_3(X,\bar a_s) 
=
 X'''\bar a_s   \hfill \cr\noalign{\smallskip}
J_4(X,\bar a_s)
=
sX^{(IV)}\bar a_s + 2X'''\bar a_s'      \hfill \cr\noalign{\smallskip}  
J_5(X,\bar a_s) 
=
s(2s-1)X^{(V)}\bar a_s + 5(2s-1)X^{(IV)}\bar a_s' + 10X'''\bar a_s''  \hfill 
\end{array}
$$

It follows that the mapping 
${\cal D}^4\to{\cal F}_4\oplus\cdots\oplus{\cal F}_0$ 
defined by (\ref{cum}) is 
$sl_2$-equivariant. Indeed, for a vector field $X \in sl_2$
(which is a polynomial in $x$ of degree $\leq2$) all the terms
$J_m(X,\bar a_s)$ in
(\ref{nac}) vanish.

Proposition 4.1 (i) is proven.

\vskip 0,3cm

Let us prove the uniqueness.
By definition, the functions $\bar a_i(x)$ are linear
expressions in $a_s(x)$ and their derivatives:
$$
\bar a_i(x)=\sum_{s,t} \alpha^s_{t}(x)a_s^{(t)}(x),
$$
where $a_s^{(t)}=d^ta_s/dx^t$, $\alpha^j_{k}(x)$ are some functions.
The fact that the normal symbol $\overline{P}_A$ is $sl_2$-equivariant
means that for a vector field $X \in sl_2$,
$\bar a_i^X=L^i_X(a_i)$.
\hfill\break
a) Substitute $X=d/dx$ to obtain that the coefficients
$\alpha^j_{k}$ does not depend on $x$;
\hfill\break
b) substitute $X=xd/dx$ to obtain the condition $j-k=i$:
$$
\bar a_i(x)=\sum_{j=4}^{i} \alpha^ja_j^{(j-i)}(x),
$$
and $\alpha^i\not=0$;
\hfill\break
c) put $\alpha^i=1$ and, finally, substitute $X=x^2d/dx$ to obtain
the coefficients from (\ref{cum}).

Proposition 4.1 (ii) is proven.

\vskip 0,3cm

The notion of normal symbol of a 4-th order differential operator
plays a central role in this paper.

\vskip 0,3cm

{\bf 4.3 Remark: the transvectants}.
Operations $J_3(X,a_s) , J_4(X,a_s) , J_5(X,a_s)$ are particular cases of 
the following remarkable bilinear operations on tensor-densities.
Consider the expressions:
$$
j_n(\phi, \psi ) 
= 
\sum_{i+j=n} (-1)^i{n \choose i}
\frac{{(2\lambda-i)! (2\mu-j)!}}{{(2\lambda-n)! (2\mu-n)!}}  
{\phi}^{(i)} {\psi}^{(j)} 
$$ 
where $\phi=\phi(x), \psi=\psi(x)$
are smooth functions.

This operations defines unique (up to a constant)
$sl_2$-equivariant mapping:
$$
{\cal F}_{\lambda} \otimes {\cal F}_{\mu} \to
 {\cal F}_{\lambda+\mu-n}
$$
Operations $j_n(\phi, \psi )$ were discovered by Gordan \cite{Gor},
they are also known as Rankin-Cohen brackets (see \cite{Ran}, \cite{Coh}).

Note, that the operations $J_m$
from the formula (\ref{nac}) are proportional to $j_n$ for $X\in{\cal F}_1,
a_s\in{\cal F}_{-s}$.

\section{Diagonalization of operator $T$}
We will obtain here an important property of
the intertwining operator (\ref{T}): in terms of normal symbol
it has a diagonal form.
We will also prove the part (i) of Theorem 1 and Theorem 4.

\vskip 0,3cm

{\bf 5.1 Proof of Theorem 1, part (i)}.
Let us define an isomorphism of modules
${\cal D}^3_{\lambda}$ and ${\cal  D}^3_{\mu}$
for $\lambda\not=0,-1,-\frac{1}{2},-\frac{1}{2}\pm \frac{\sqrt 21}{6}$. 
Associate to $A \in {\cal D}^3_{\lambda}$ the operator $T(A) \in 
{\cal  D}^3_{\mu}$:
$$
T:\;a_3\frac{d^3}{dx^3}+a_2\frac{d^2}{dx^2}+a_1\frac{d}{dx}+a_0\;\longmapsto\;
a_3^T\frac{d^3}{dx^3}+a_2^T\frac{d^2}{dx^2}+
a_1^T\frac{d}{dx}+ a_0^T
$$
 such that its standard symbol 
$$
{\overline P} _{T(A)} = 
\xi^3 {\overline {a^T}_3(x)}+
\xi^2{\overline {a^T}_2(x)}+
\xi{\overline {a^T}_1(x)}+
 {\overline {a^T}_0(x)}
$$ 
is given by:
\begin{equation}
\matrix{
{\overline {a^T}_3(x)}
=
\bar a_3(x) \hfill\cr \noalign{\smallskip}
{\overline {a^T}_2(x)}
= 
\displaystyle\frac {{2\mu +1}}{{2\lambda +1}} \bar a_2(x) \hfill\cr  \noalign{\smallskip}
{\overline {a^T}_1(x)}
=
\displaystyle\frac {{3{\mu}^2+3\mu -1}}{{3{\lambda}^2+3\lambda -1}} \bar a_1(x) \hfill\cr  \noalign{\smallskip}
{\overline {a^T}_0(x)}
=
\displaystyle\frac{{{\mu}(\mu +1)(2\mu +1)}}{{{\lambda}(\lambda +1)(2\lambda +1)}} \bar a_0(x)
\hfill\cr  \noalign{\smallskip}
}
\label{TTT}
\end{equation}
It follows immediately from the formula (\ref{nac}), 
that this formula defines an isomorphism 
of $\Vect(\bf R)$-modules:
$
T:{\cal D}^3_{\lambda}\cong{\cal D}^3_{\mu}.
$

Theorem 1 (i) is proven.

\vskip 0,3cm

{\bf 5.2 Proof of Theorem 4}.
Let us show that the operator (\ref{TTT}) in terms of symmetric expressions
of Lie derivatives is given by
(\ref{int}).

The first equality in (\ref{int}) coincides with the last equality in (\ref{TTT}).

1) Consider a first order operator of a Lie derivative 
$$
L_X^{\lambda}=X(x)\frac{d}{dx}-\lambda X^{\prime}(x).
$$
Its normal symbol defined by (\ref{cum}) is
$$
\overline{P}_{L_X^{\lambda}}=\xi X(x).
$$
One obtains the second equality of the formula (\ref{int}).

2) The anti-commutator
$$
[L_X^{\lambda},L_Y^{\lambda}]_+=2XY\frac{d^2}{dx^2}+(1-2\lambda)(XY)^{\prime}\frac{d}{dx}
-\lambda(XY^{\prime\prime}+X^{\prime\prime}Y)+2\lambda^2X^{\prime}Y^{\prime}
$$
has the following normal symbol:
$$
\overline{P}_{[L_X^{\lambda},L_Y^{\lambda}]_+}=
2\xi^2XY-
\frac{2}{3}\lambda(\lambda+1)(XY^{\prime\prime}+X^{\prime\prime}Y-X^{\prime}Y^{\prime})
$$
which also following from (\ref{cum}).
The third equality of (\ref{int}) follows now from the second and the fourth ones
of (\ref{TTT}).

3) The normal symbol of a third order expression 
$[L_X^{\lambda},L_Y^{\lambda},L_Z^{\lambda}]_+:=\break\hbox{Sym}_{X,Y,Z}(L_X^{\lambda}L_Y^{\lambda}L_Z^{\lambda})$
can be also easily calculated from (\ref{cum}). The result is:
$$
\matrix{
\overline{P}_{[L_X^{\lambda},L_Y^{\lambda},L_Z^{\lambda}]_+}=6\xi^3 XYZ \hfill\cr  \noalign{\smallskip}
\;\;\;\;\;\;\;\;\;\;\;\;\;\;\;\;\;\;\;\;
-(3\lambda^2+3\lambda-1)\xi 
(XYZ^{\prime\prime}+XY^{\prime\prime}Z+X^{\prime\prime}YZ-
\frac{1}{5}(XYZ)^{\prime\prime})\hfill\cr  \noalign{\smallskip}
\;\;\;\;\;\;\;\;\;\;\;\;\;\;\;\;\;\;\;\;
-\lambda(\lambda+1)(2\lambda+1)
(XYZ^{\prime\prime\prime}+XY^{\prime\prime\prime}Z+X^{\prime\prime\prime}YZ)\hfill\cr
}
$$
This formula implies the last equality of (\ref{int}).

\vskip 0,3cm

{\bf 6.3 Remarks}.  a) The normal symbols of $[L_X^{\lambda},L_Y^{\lambda}]_+$
and $[L_X^{\lambda},L_Y^{\lambda},L_Z^{\lambda}]_+$ are given by
very simple and harmonic expressions
(which implies the diagonal form (\ref{int}) of operator $T$).
It would be interesting to understand better the geometric reason of this fact.

b)
Comparing the formul{\ae} (\ref{int}) and (\ref{TTT}), one finds a
coincidence between coefficients. This fact shows that, in some sense,
the symmetric expressions of Lie derivatives and the normal symbol 
represent the same thing in terms of differential operators and in terms
of polynomial functions on $T^*{\bf R}$, respectively. We do not see any reason {\it a-priori}
for this remarkable coincidence.

\section{Uniqueness of operator $T$}
In this section we prove that the isomorphism $T$ 
defined by the formula (\ref{TTT}) is unique (up to a constant).
We also show that in the higher order case $k\geq4$ there is no analogues
of this operator.

\vskip 0,3cm

{\bf 6.1 Proof of Theorem 3}.
The normal symbol of an operator $A\in{\cal D}^3_{\lambda}$ {\it is}
at the same time a normal symbol of $T(A)\in{\cal D}^3_{\mu}$,
since operator $T$ is equivariant.
The normal symbol is unique up to normalization (cf. Proposition 4.1, part (ii)),
therefore
the polynomial ${\overline P} _{T(A)}(x,\xi)$
defined by the formula (\ref{cum}), is of the form:
$$
{\overline P} _{T(A)}(x,\xi)=
\alpha_3\xi^3\bar a_3(x)+\alpha_2\xi^2\bar a_2(x)+\alpha_1\xi\bar a_1(x)+\alpha_0\bar a_0(x),
$$
where $\alpha_i\in \bf R$ are some constants
depending on $\lambda$ and $\mu$.
Choose $\alpha_3=1$. It follows immediately from the formula
(\ref{nac}) (after substitution $a_4\equiv0$)
that the formula (\ref{TTT}) gives the unique choice of the constants
$\alpha_2,\alpha_1,\alpha_0$ such that
operator $T$ is equivariant.

Theorem 3 is proven.

\vskip 0,3cm

{\bf 6.2 Proof of Theorem 2}.
Suppose now that $\Phi : {\cal D}^4_{\lambda} \rightarrow {\cal  D}^4_{\mu}$ 
is an isomorphism. In the same way, it follows that in terms of normal symbols,
operator $\Phi$ 
is diagonal.
More precisely, if $A \in {\cal D}^4_{\lambda}$, then the normal symbol 
of the operator $\Phi (A) \in {\cal  D}^4_{\mu}$ is: 
$$
\overline{P}_{\Phi (A)}(x,\xi)
= 
\alpha _4\xi^4{\bar a_4}(x)+
\alpha_3\xi^3\bar a_3(x)+\alpha_2\xi^2\bar a_2(x)+\alpha_1\xi\bar a_1(x)+\alpha _0{\bar a_0}(x),
$$
where $\bar a_i$ are the components of the normal symbol of $A$, $\alpha
_i \in \bf R$. 
The condition of equivariance implies : 
$$
\matrix{
\displaystyle\frac{{\alpha _2}}{{\alpha _0}}
=
\frac{{{\lambda}(\lambda +1)}}{{{\mu}(\mu +1)}},\hfill&
\displaystyle\frac{{\alpha _3}}{{\alpha _0}}
=
\frac{{{\lambda}(\lambda +1)(2\lambda +1)}}{{{\mu}(\mu +1)(2\mu +1)}}
\hfill\cr\noalign{\smallskip}
\displaystyle\frac{{\alpha _4}}{{\alpha _0}}
=
\frac{{{\lambda}(\lambda +1)(12{\lambda}^2+12\lambda +11)}}
{{{\mu}(\mu +1)(12{\mu}^2+12\mu +11)}},\hfill&
\displaystyle\frac{{\alpha _4}}{{\alpha _2}} 
= 
\frac {{6{\lambda}^2+6\lambda -5}}{{6{\mu}^2+6\mu -5}}\hfill\cr\noalign{\smallskip}
\displaystyle\frac{{\alpha _3}}{{\alpha _1}}
=
\frac {{3{\lambda}^2+3\lambda -1}}{{3{\mu}^2+3\mu -1}},\hfill&
\displaystyle\frac{{\alpha _4}}{{\alpha _1}}
=
\frac {{{\lambda}(\lambda +1)(2\lambda +1)}}
{{{\mu}(\mu +1)(2\mu +1)}}\hfill\cr\noalign{\smallskip}
} 
$$

This system of equation has solutions if and only if $\lambda = \mu$ or $\lambda+\mu 
= -1$. For $\lambda = \mu$ one has:
$\alpha_0=\alpha_1=\alpha_2=\alpha_3=\alpha_4$ and for $\lambda+\mu 
= -1$ one has: $\alpha_0=-\alpha_1=\alpha_2=-\alpha_3=\alpha_4$.

Theorem 2 is proven for $k=4$.

\vskip 0,3cm 

Theorem 2 follows now from one of the results of \cite{lec}:
given  an isomorphism $\Phi:{\cal D}^k_{\lambda}\to {\cal D}^k_{\mu}$,
then the restriction of $\Phi$ to ${\cal D}^4_{\lambda}$ is an 
isomorphism of $\Vect(\bf R)$-modules:
${\cal D}^4_{\lambda}\to{\cal D}^4_{\mu}$.
(To prove this, it is sufficient to suppose equivariance of
$\Phi$ with respect to the affine algebra with generators $\frac{d}{dx}, x\frac{d}{dx}$,
see \cite{lec}).

This implies that $\lambda = \mu$ or $\lambda+\mu = 
-1$. 

Theorem 2 is proven.

\section{Relation with the cohomology group\break
$H^1 (\Vect(\bf R);\hbox{Hom}({\cal F}_{\lambda},{\cal F}_{\mu}))$}

  The problem of isomorphism 
of $\Vect({\bf R})$-modules ${\cal D}^k_{\lambda}$
for different values of $\lambda$
 is related to the 
first cohomology group
 $H^1(\Vect({\bf R});\hbox{Hom}({\cal F}_{\lambda},{\cal F}_{\mu}))$. 
This cohomology group has already been
calculated by B.L. Feigin and D.B. Fuchs (in the case of formal series) \cite{fei}.

\vskip 0,3cm

{\bf 7.1 Nontrivial cocycles}.The relation of $\Vect(\bf R)$-action on the space of differential 
operators and the cohomology groups $H^1 (\Vect(\bf R);\Hom ({\cal F}_{\lambda}, {\cal 
F}_{\mu}))$ is given by the following construction. 

Let us associate to the bilinear mappings $J_m$
defined by the formula (\ref{nac}),
a linear mapping 
$
c_m : \Vect(\bf R) \rightarrow \Hom ({\cal F}_s, {\cal F}_{s+1-m})   
$: 
$$
c_m(X)(a) := J_m(X,a), 
$$ 
where $a \in {\cal F}_s$.

A remarkable property of transvectants $J_3$ and $J_4$ is:

\proclaim Lemma 7.1. For each value of $s$,
the mappings $c_3$ and $c_4$ are non-trivial cocycles on $\Vect(\bf R)$:
\hfill\break
(i) $
c_3\in Z^1(\Vect(\bf R);\Hom ({\cal F}_s, {\cal F}_{s-2})),
$\hfill\break
(ii) $
 \;c_4\in Z^1(\Vect(\bf R);\Hom ({\cal F}_s, {\cal F}_{s-3})).
$\par

{\bf Proof}. From the fact that the formula (\ref{nac}) defines a $\Vect(\bf R)$ 
action one checks that for any $X ,Y \in \Vect(\bf R)$ : 
$$ 
[L_X,c_m(Y)] - [L_Y,c_m(X)] = c_m([X,Y])   
$$ 
with $m = 3 , 4$.
This means, that $c_3$ and $c_4$ are cocycles.

The cohomology classes $[c_3], [c_4]\not=0$. 
Indeed,
verify that $c_3$ and $c_4$ are 
cohomological to the non-trivial cocycles: 
$$
\begin{array}{rcl} 
\widetilde c_3(X)(a) 
&=&
X'''a + 2X''a' , \\ \noalign{\smallskip}
\widetilde c_4(X)(a)  
&=& 
X'''a' + X''a'' \\ \noalign{\smallskip} 
\end{array}
$$
 from \cite{fei}.

Lemma 7.1 is proven.

\vskip 0,3cm

{\bf 7.2 Proof of Theorem 1, part (ii)}. 
First, remark that for every $\alpha _1, \alpha _2, \alpha _3 \in \bf R$, the formula
$$
\matrix{ 
\rho_X(a_3)  
&=& 
L_X^3(a_3)  \hfill  \cr\noalign{\smallskip}
\rho_X(a_2)  
&=& 
L_X^2(a_2) \hfill \cr\noalign{\smallskip}
\rho_X(a_1)  
&=& 
L_X^1(a_1) + {\alpha _1}J_3(X,a_3)     \hfill  \cr\noalign{\smallskip}
\rho_X(a_0)  
&=& 
L_X^0(a_0) + {\alpha _2}J_3(X,a_2) + {\alpha _3}J_4(X,\bar a_3)   \hfill 
}
$$
defines a $\Vect(\bf R)$-action. 
Indeed, this formula coincides with (\ref{nac})
in the case $a_4\equiv0$ and for the special values of
$\alpha _1, \alpha _2, \alpha _3$, however, the constants
$\alpha _1, \alpha _2, \alpha _3$
are independent.

The $\Vect(\bf R)$-action $\rho$ is a {\it non-trivial 3-parameter deformation}
of the standard action on the direct sum 
${\cal F}_3\oplus{\cal F}_2\oplus{\cal F}_1\oplus{\cal F}_0$.

The fact, that
the cocycles $c_3$ and $c_4$ are non-trivial,
is {\it equivalent} to the fact that
the defined $\Vect(\bf R)$-modules with:

1) $\alpha_1,\alpha_2,\alpha_3 \not =0$, 

2) $\alpha_1=0,\alpha_2\not =0,\alpha_3\not =0$,

3) $\alpha_1\not =0,\alpha_2=0,\alpha_3\not =0$,

4) $\alpha_1\not =0,\alpha_2 \not=0,\alpha_3=0$,
\hfill\break
are not isomorphic to each other.

The $\Vect(\bf R)$-modules ${\cal D}^3_{\lambda}$
(given by the formula (\ref{nac}) with $a_4\equiv0$)
corresponds to the case 1) for general values of $\lambda$,
to the case 2) for $\lambda=-{1\over 2}\pm \frac{\sqrt 21}{6}$,
to the case 3) for $\lambda=-{1\over 2}$
and to the case 4) for $\lambda=0,-1$.
Therefore,
one obtains 5 critical values of the degree for which 
$\Vect(\bf R)$-module structure on the space of third order operators is special.

Theorem 1 (ii) is proven. 

\vskip 0,3cm

{\bf Remark}. For each value of $\lambda$
at least one of constants $\alpha_1,\alpha_2,\alpha_3 \not =0$.
This implies that the module ${\cal D}^3_{\lambda}$
is not isomorphic to the direct sum
${\cal F}_3\oplus{\cal F}_2\oplus{\cal F}_1\oplus{\cal F}_0$.

\section{Explicit formula for the intertwining operator}

We give here the explicit formula for the operator (\ref{T})
intertwining $\Vect(\bf R)$-actions on ${\cal D}^3$
which follows from the expression for the operator $T$
in terms of the normal symbol
(\ref{TTT}).

For every $A \in {\cal D}^3_{\lambda}$, the operator $T(A) \in 
{\cal  D}^3_{\mu}$
$$
T(A)=
a_3^T\frac{d^3}{dx^3}+a_2^T\frac{d^2}{dx^2}+
a_1^T\frac{d}{dx}+ a_0^T
$$
is given by the following formula:

\begin{equation}
\matrix{
a^T_3=a_3 \hfill\cr\noalign{\bigskip}
a^T_2=
\displaystyle\frac{2\mu+1}{2\lambda+1}a_2+
\frac{3(\mu-\lambda)}{2\lambda+1}a_3^{\prime} \hfill\cr\noalign{\bigskip}
a^T_1=
\displaystyle\frac{3\mu^2+3\mu-1}{3\lambda^2+3\lambda-1}a_1 \hfill\cr\noalign{\smallskip}
\;\;\;\;\;\;+
\displaystyle\frac{(\lambda-\mu)(\mu(12\lambda-1)-\lambda+3)}
{2(2\lambda+1)(3\lambda^2+3\lambda-1)}a_2^{\prime} \hfill\cr\noalign{\smallskip}
\;\;\;\;\;\;+
\displaystyle\frac{3}{2}\frac{\mu^2(5\lambda-1)-\mu(6\lambda^2+\lambda-1)+
\lambda^3+2\lambda^2-\lambda}
{(2\lambda+1)(3\lambda^2+3\lambda-1)}a_3^{\prime\prime} \hfill\cr\noalign{\bigskip}
a^T_0
=
\displaystyle\frac{\mu(\mu+1)(2\mu+1)}{\lambda(\lambda+1)(2\lambda+1)}
a_0\hfill\cr\noalign{\smallskip}
\;\;\;\;\;\;-
\displaystyle\frac{\mu^3(3\lambda+5)-\mu^2(3\lambda^2-6)-\mu(5\lambda^2+6\lambda)}
{(\lambda+1)(2\lambda+1)(3\lambda^2+3\lambda-1)}a_1^{\prime} \hfill\cr\noalign{\smallskip}
\;\;\;\;\;\;+
\displaystyle\frac{\mu^3(3-\lambda)-\mu^2(6\lambda^2+7\lambda-5)+
\mu(7\lambda^3+4\lambda^2-5\lambda)}
{2(\lambda+1)(2\lambda+1)(3\lambda^2+3\lambda-1)}a_2^{\prime\prime} \hfill\cr\noalign{\smallskip}
\;\;\;\;\;\;-
\displaystyle\frac{\mu^3(3\lambda^2+1)-3\mu^2(\lambda^2+2\lambda-1)
-\mu(3\lambda^4-3\lambda^3-5\lambda^2+3\lambda)}
{2(\lambda+1)(2\lambda+1)(3\lambda^2+3\lambda-1)}a_3^{\prime\prime\prime}\hfill\cr\noalign{\bigskip}
}
\label{exp}
\end{equation}

We de not prove this formula
since we do not use it in this paper.

\vskip 0,3cm

{\bf Remarks}. 1) If $\lambda=\mu$, then
the operator $T$ defined by this formula is identity,
if $\lambda+\mu=-1$, then $T$ is the operator of conjugation.

2) The fact that operator $T$ is equivariant implies that
{\it the formula (\ref{exp}) does not depend on the choice of the coordinate $x$}.

\section{Discussion and final remarks}

Let us give here few examples
and applications of the normal symbols (\ref{cum}).

\vskip 0,3cm

{\bf 9.1 Examples}. The notion of normal symbol 
 was introduced in \cite{duv}
in the case of second order differential operators.
In this case, for $\lambda=\frac{1}{2}$ (operators on $-\frac{1}{2}$-densities), 
the normal symbol
(\ref{cum}) is just the standard total symbol:
$\bar a_2=a_2,\bar a_1=a_1,\bar a_0=a_0$.
This corresponds to the classical example of second order operators
on $-\frac{1}{2}$-densities (cf. the footnote at the introduction).

In the same way, the semi-integer values $\lambda=1,\frac{3}{2},2,\dots$
corresponds to particularly simple expressions of the normal symbol
for operators of order $k=3,4,5,\dots$.

\vskip 0,3cm

{\bf 9.2 Modules of second order differential operators on $\bf R$}. 
a) The module of second order operators ${\cal D}^2_{\lambda}$
with $\lambda=0,-1$,
is decomposed to a direct sum:
$$
{\cal D}^2_0\cong{\cal D}^2_{-1}\cong{\cal F}_2\oplus{\cal F}_1\oplus{\cal F}_0.
$$
Indeed, the coefficients $\bar a_2,\bar a_1,\bar a_0$
transform as tensor-densities (cf formula (\ref{nac})).
This module is special: ${\cal D}^2_{\lambda}$ is not isomorphic
to ${\cal D}^2_0$ for $\lambda\not=0,-1$ (see \cite{duv}).

b) For every $\lambda, \mu \not=0,-1$, ${\cal D}^2_{\lambda}\cong{\cal D}^2_{\mu}$.

\vskip 0,3cm

{\bf 9.3 Operators on $\frac{1}{2}$-densities}.
For every $k\geq3$, the module ${\cal D}^k_{-{1\over2}}$ (corresponding to $\frac{1}{2}$-densities)
is special. It is decomposed into
a sum of submodules: of symmetric operators and of skew-symmetric operators.

\vskip 0,3cm

{\bf 9.4 Normal symbol and Weil symbol}.
The Weil quantization defines a 1-parameter family of mappings from
the space of polynomials $C[\xi,x]$ to the space of differential
operators on $\bf R$ with polynomial coefficients.
One associate to a polynomial the symmetric expression in
$\hbar\frac{d}{dx}$ and $x$:
$F(\xi,x)\mapsto
\hbox{Sym}F(\hbar\frac{d}{dx},x)$.
This (1,1)-correspondence between differential operators and
polynomials is $sl_2$-equivariant.
However, in the Weil quantization
the action of the Lie algebra $sl_2$ on differential operators 
is generated by $x^2,x\frac{d}{dx}+\frac{d}{dx}x,\frac{d^2}{dx^2}$
and therefore, is completely different from the normal symbol.

\vskip 0,3cm

{\bf 9.5 Automorphic (pseudo)differential operators}.
The notion of canonical symbol is related 
(and in some sense inverse) to the construction of
the recent work
P. Cohen, Yu. Manin and D. Zagier \cite{cmz}
of a $PSL_2$-equivariant (pseudo)differential operator
associated to a holomorphic tensor-density on the upper half-plane.

\vskip 0,3cm

{\bf 9.6 Exotic $\star$-product}.
Another way to understand this $sl_2$-equivariant 
correspondence between linear differential operators and polynomials in $\xi,x$
leads to
a $\star$-product on the algebra of Laurent polynomials on $T^*{\bf R}$
which is
not equivalent to the standard Moyal-Weil quantization
(see \cite{qua}).

\vskip 0,3cm

{\it Acknowledgments}. It is a pleasure to acknowledge enlightening
discussions with C. Duval,
P.B.A. Lecomte, Yu.I. Manin and E. Mourre. 

\vskip 0,5cm


\end{document}